\documentclass[12pt]{article}
\usepackage{amsmath,amsfonts,amssymb,amsthm}

\def\noi{\noindent}

\def\pf{\noi{\bf Proof.\ \,}}
\def\eop{{$\square$}}
\def\qed{{$\square$}}

\def\labtt#1{\label {#1}  }
\def\labttr#1{\label {#1} \rm }

\def\refpp#1{(\ref {#1})}

\def\a{\alpha}
\def\b{\beta}

\def\l{\lambda}

\def\vep{\varepsilon}

\def\L{\Lambda}

\def\CC{{\mathbb C}}
\def\FF{{\mathbb F}}
\def\MM{{\mathbb M}}
\def\QQ{{\mathbb Q}}
\def\RR{{\mathbb R}}

\def\TT{{\mathbb T}}
\def\ZZ{{\mathbb Z}}

\def\la{\langle}
\def\ra{\rangle}
\def\<{\langle}
\def\>{\rangle}

\def\l{{\lambda}}

\def\half{{1 \over 2}}

\def\dual#1{#1^*}        

\def\dg#1{{\cal D}({#1})}

\begin{document}

\newtheorem{thm}{Theorem}[section]
\newtheorem{prop}[thm]{Proposition}
\newtheorem{lem}[thm]{Lemma}
\newtheorem{rem}[thm]{Remark}
\newtheorem{coro}[thm]{Corollary}
\newtheorem{conj}[thm]{Conjecture}
\newtheorem{de}[thm]{Definition}
\newtheorem{hyp}[thm]{Hypothesis}

\newtheorem{nota}[thm]{Notation}
\newtheorem{ex}[thm]{Example}
\newtheorem{proc}[thm]{Procedure}



%

\begin{center}\end{center}

\begin{center}
{\Large \bf
A moonshine path for $5A$  and associated lattices of ranks 8 and 16}
\begin{center}

\vspace{10mm}
Robert L.~Griess Jr.
\\[0pt]
Department of Mathematics\\[0pt] University of Michigan\\[0pt]
Ann Arbor, MI 48109  USA  \\[0pt]
{\tt rlg@umich.edu}\\[0pt]
\vskip 1cm

Ching Hung Lam
\\[0pt]
Institute of Mathematics \\[0pt]
Academia Sinica\\[0pt]
Taipei 10617, Taiwan\\[0pt]
{\tt chlam@math.sinica.edu.tw}\\[0pt]
\vskip 1cm
\end{center}


\end{center}

\begin{abstract}  
We continue the program, begun in \cite{gl3cpath}, to make a moonshine path between a node of the extended $E_8$-diagram and the Monster simple group.
Our goal is to  provide a context for observations of McKay, Glauberman and Norton
by
realizing their theories in a more concrete form.
In this article, we treat the $5A$-node.  Most work in this article is a study of certain lattices of ranks 8, 16 and 24 whose determinants are powers of 5.
\end{abstract}

\newpage

\tableofcontents

\newpage

\section{Introduction}

 In \cite{gl3cpath},  we  initated a  program for establishing
moonshine paths between the extended $E_8$-diagram  and the monster simple group and treated the case of the $3C$-node in detail.
The present article treats the case of the $5A$-node, which corresponds to
 a pair of Monster involutions in class $2A$ which generate a dihedral group of order 10.

Such a group corresponds naturally to a
dihedral group of order 10 which acts on the rank 24 Leech lattice
and which is generated by pairs of trace 8 involutions.
Most work in this article is a study of
lattices, in particular certain rootless rank 8 and 16 lattices of determinant $5^4$  and their embeddings in Niemeier lattices and other related embeddings.  We give attention to the role of $EE_8$-sublattices  because of their occurrence within $Q$ and $\L$ and their role in VOA theory \cite{glee8}.

The introduction of  \cite{gl3cpath} has a detailed discussion of context,
which involves lattices, VOAs, Lie theory and finite groups.  Our desire is to
take mystery out of McKay's $E_8$-Monster connection.  In particular, we note
the Glauberman-Norton theory \cite{gn} with its interesting observation that
particular centralizers in the Monster involve ``half'' the Weyl group of a
node in the extended $E_8$-diagram.  In \cite{gl3cpath}, we gave an
interpretation within our moonshine path theory for why just half a Weyl group
occurrs in the $3C$-case.
The present
article treats the $5A$-case.
Work on other nodes is planned.

\begin{thm}\labtt{mainth1} Suppose that
$R$ is a rank 8 rootless even lattice such that $\dg R \cong 5^4$.   Then

(i) $R$ is unique up to isometry.

(ii)  $O(R)\cong O^+(4,5)$.

(iii) $O(R)$ acts faithfully on $\dg R$ and $O(R)$ contains no SSD involutions.  

\end{thm}

\begin{thm}\labtt{mainth2} Suppose that
$Q$ is a rank 16 rootless even lattice such that $\dg Q \cong 5^4$.  Then

(i) $Q$ is unique up to isometry.

(ii) There exists an embedding $O(Q)\rightarrow Frob(20) \times O^+(4,5)$ such
that the image has index 2 and contains neither direct factor..

(iii) The action of $O(Q)$ on $\dg Q$ induces $O(\dg Q)\cong O^+(4.5)$ and the
kernel is isomorphic to $Dih_{10}$.   The set of SSD-involutions in $O(Q)$
generates this kernel.  
\end{thm}


\begin{thm}\labtt{mainth3}
Let $\L$ be a Leech lattice.  Let $R, Q$ be integer lattices of respective ranks 8, 16, such that
$\dg R \cong \dg Q\cong 5^4$.

(i)  In $\L$, the set of sublattices isometric to $R$ forms an orbit under $O(\L )$.

(ii) In $\L$, the set of sublattices isometric to $Q$ forms an orbit under $O(\L )$.
\end{thm}


%


\begin{coro}\labtt{mainth4}
Let $D$ be the dihedral group generated by a pair of defining $EE_8$-involutions in
$Q\cong DIH_{10}(16)$. Then $N:=N_{O(\L )}(D)$ has shape $N/D \cong
O^+(4,5)$ and $N$ is isomorphic to the unique subgroup of index 2 in $Frob(20)
\times O^+(4,5)$ which contains the subgroup $Dih_{10} \times GL(2,5)$ and
does not contain either direct factor of $Frob(20) \times O^+(4,5)$. The action of
$N$ on $Q$ is faithful and the action of $N$ on $R$ has kernel $D$.
\end{coro}


\begin{coro}\labtt{maincorossd}
In $O(R)$, the set of SSD involutions is empty.  In O(Q), the set of SSD involutions is the set of five $EE_8$-involutions in the normal dihedral group.
\end{coro}
\pf
Use \refpp{mainth1}, \refpp{mainth2} and \refpp{trivialaction}.
\eop


\begin{coro}\labtt{maincoro1} The centralizer in $O(\L)$ of $D$ is isomorphic to
$$(SL(2,5){\circ}SL(2,5)){:}2,$$ where the outer elements of order 2 normalize each of the two quasisimple central factors under conjugation and induce noninner automorphisms.
\end{coro}


\medskip

%


\begin{thm}\labtt{mainth5}  Let $R, Q$ be as in \refpp{mainth1},\refpp{mainth2}.  There are three orbits of $O(Q\perp R)$ on the set of Niemeier overlattices of $Q\perp R$.
Let $L$ be such an overlattice and define the integer $e$ by $5^e:=|\QQ \otimes R \cap L:R|$.  The three orbits are distinguished by the common value of $e\in \{0, 1, 2\}$ for their members.  We have
$L\cong \L, {\mathcal N}(A_4^6), E_8^3$ for $e=0,1,2$, respectively.
The stabilizer in $O(Q)\times  O(R)$ of
a Niemeier overlattice is a subgroup $H$ of the following form (we identify $O(Q)\times  O(R)$ with $O(Q\perp R)$ in the obvious way):

$e=0$:  $H\cap O(R)=1$, $H\cap O(Q)\cong Dih_{10}$, $H \cong O(Q)$ (so $H/(H\cap O(Q))\cong O^+(4,5)$;

$e=1$: $H\cap O(R)=5$, $H\cap O(Q)\cong Dih_{10} \times 5$, $H/(H\cap O(Q))\cong O^+(2,2)\cong 4 \wr 2$.

$e=2$: $H=(H\cap O(Q)) \times (H\cap O(R))$, $H \cap O(Q)\cong Dih_{10}.5{:}GL(2,5)$, $H\cap O(R)\cong 5{:}GL(2,5)$.
\end{thm}

\medskip

We obtain a structure analysis of $C_{O(\L)}(D)$
 and in particular show that $C_{O(\L)}(D)/O_2(C_{O(\L)}(D))$ is isomorphic to ``half'' of the Weyl group of
type $A_4^2$.
These results
may be lifted to a statements \cite{gn} about the Monster as
follows. There exist suitable $x, y, z$ in the Monster ($x, y$ in $2A$, $xy$ in
$5A$, $z\in 2B \cap C_{\MM}(x,y)$) so that $C_{\MM}(x,y,z)$ has a homomorphism onto
$C_{O(\L)}(D)/O_2(C_{O(\L)}(D))$.
In analogy with \cite{gl3cpath}, we point out a presence of triality.

As in \cite{gl3cpath}, a study of associated conformal vectors of central
charge $\half$  leads to embedding results for certain VOAs. We shall refer to
\cite{gl3cpath}, Section 2   and Appendix B,   for details.




%







\subsection{Acknowledgments}

The first author thanks the United States National Science Foundation (NSF DMS-0600854)
and National Security Agency (NSA H98230-10-1-0201)
for financial support
and the Academia Sinica for hospitality during a visit, December, 2009.

The second author thanks  National Science Council (NSC 97-2115-M-006-015-MY3)
and National Center for Theoretical Sciences of Taiwan for financial support
and University of Michigan for hospitality during a visit, March, 2010.

We thank Gabriele Nebe, Rudolf Scharlau and Rainer Schulze-Pillot for helpful
advice about integer lattices. After we obtained our results, we learned that
the characterization of $R$ as in \refpp{mainth1} had been known
\cite{scharlauhemkemeier}. Also, see the earlier articles \cite{costellohsia},
\cite{maass}.

Quite recently, a study of the genus of rank 16, discriminant group $5^4$ integral
lattices was completed by Scharlau and Hemkemeier, using the
software package ``tn'' \cite{tn}.   Those authors found that the genus has 848 members,
just one of which is rootless.  This is in agreement with our \refpp{mainth2}.
 A reference to this lattice and its isometry group  is row 11 in  \cite{pn}, page 86 and \cite{rsrsp}.

\newpage

\begin{center}
{\bf \large Notation and Terminology}

\vspace{0.3cm}
\begin{tabular}{|c|c|c|}
   \hline
\bf{Notation}& \bf{Explanation} & \bf{Examples in text}  \cr
   \hline\hline
$2A, 2B, 3A, \dots$ & conjugacy classes of the Monster: &  Introduction\cr
                  & the first number denotes the order & \cr
                  & of the elements and the second letter & \cr
                  & is arranged in descending order of & \cr
                  & the size of the centralizers &\cr \hline
$A_1,  \cdots , E_8$ &root lattice for root system &  Sec. 2\cr
& $\Phi_{A_1}, \dots ,
\Phi_{E_8}$& \cr
   \hline
$AA_1, \cdots ,$ & lattice isometric to $\sqrt 2$ times &\cr
$ EE_8$ &the lattice  $A_1, \cdots,  E_8$ & Page 2\cr
\hline
$A\circ B$ & central product of  groups $A$ and $B$ & Coro \ref{fourdimog }     \cr \hline
$\dg L$ & discriminant group of integral  & Notation \ref{dg}\cr
 &  lattice $L$: $\dual L/L$& Lemma \ref{maximalWitt}\cr      \hline
 $EE_8$-involution & SSD involution whose negated space &  \refpp{q4} \cr
& is isometric to $EE_8$ & \cr \hline
 $G', G'' ,\dots $ & first, second, \dots commutator subgroup  &   \cr
& of the group $G$ &  \refpp{centdniema4to6} \cr \hline
 $G^{(\infty )}$ & the terminal member of the derived  &   \cr
& series (commutator series) of the group $G$ &  \refpp{centdniema4to6} \cr \hline
$h$ & a fixed point free automorphism& Remark \ref{r=(h-1)E8}\cr
& of $E_8$ of order $5$ & \cr \hline
 $L(a,\cdots)$ & a lattice created by a gluing of    & \refpp{q7}\cr
 & 2-spaces from orbit of 2-spaces &\cr
 & from orbit of type $Orbit(a, \cdots)$ & \cr\hline
  Niemeier  & a rank 24 even unimodular lattice
  & Lemma \ref{el3}  \cr
  lattice&&\cr \hline
 $\mathcal{N}(X)$ &
Niemeier lattice whose root system& Lemma \ref{el3}\cr
 &has type $X$&               \cr \hline
  $O(X)$ & the isometry group of  & \refpp{r13}, \cr
  &the quadratic space $X$&  \refpp{q7.5} \cr\hline
  $Orbit(a,\cdots)$ & orbit of $O^+(4,q)$ on 2-spaces:  & \refpp{q6}\cr
&  $a$ denotes  the dimension of & \cr
&the radical of the 2-space &\cr \hline
root group & subgroup of group of Lie type  & \refpp{centrootgroup}\cr
& associated to a root (see \cite{Carter})  &\cr \hline
 \end{tabular}
\end{center}

%




\section{A preliminary observation}
%
\begin{nota}\label{dg}
Let $L$ be an integral lattice. We shall denote the discriminant group $L^*/L$
 of $L$ by $\dg L$.
\end{nota}

\begin{lem}\labtt{maximalWitt}
Let $R, Q$ be integer lattices of respective ranks 8, 16, such that
$\dg R \cong \dg Q\cong 5^4$.
(i)
The quadratic spaces $\dg R$ and $\dg Q$ have maximal Witt index.

(ii)
In each case there exist overlattices of
the same rank which are even and unimodular.
\end{lem}
\pf
The assumptions in \refpp{mainth1} imply that the quadratic space $\dg R$ has
maximal Witt index.  The reason is that there exists an $R$ with all these
properties (see \refpp{r8}, for example) and all lattices with these properties
have the same genus. See \cite{splag}, p.386, and \cite{scharlauhemkemeier}. A
similar discussion applies to the lattices of \refpp{mainth2}; see \cite{glee8} or
\refpp{qine8e8.2} or just observe that $Q$ and $E_8 \perp R$ are in the same genus.
\eop

\section{Uniqueness for $R$}

\begin{nota}\labttr{r6}
We let $R$ be a rootless rank 8  even integral lattice.
\end{nota}

In this section, we shall establish the uniqueness for $R$ and prove Theorem
\ref{mainth1}.

%

\begin{lem}\labtt{r8}
Let $R$ be as in \refpp{r6}.

(i) Then $R$ is contained in a copy of the $E_8$ lattice in $\QQ\otimes R$.

(ii) The orthogonal group on $\QQ\otimes R$ has one orbit on the set $\mathcal
A$ of pairs $(X,Y)$ where $X\cong R, Y\cong E_8$ and $X\le Y$. In particular, the
conditions of \refpp{r6} characterize $R$ up  to isometry.

(iii) If $X\cong R$, the number of $Y\cong E_8$ so that $(X,Y)\in \mathcal A$ is
12.

(iv) If $Y\cong E_8$, the number of $X\cong R$ so that $(X,Y)\in \mathcal A$ is
$\frac {|Weyl(E_8)|}{|5 \times GL(2,5)|} = 2^{9}3^4{\cdot 7}$.
\end{lem}
\pf (i) Since $\dg R$ has Witt index 2, there exists an integral lattice $J$ so that
$J\ge R$ and $|J:R|=25$.  Since $R$ is even and its index is odd, $J$ is even.  By
the well-known characterization of $E_8$, $J\cong E_8$. (See \cite{gre8} for a
survey of characterizations of $E_8$.)

(ii) The second statement follows from \refpp{r5}.

(iii) The number  of totally singular subspaces in $\dg R$ is 12. A
lattice between $R$ and $\dual R$ which is even and unimodular corresponds to
such a subspace.   The lattice between $R$ and $\dual R$ corresponding to such
a subspace is even and unimodular, hence isometric to $E_8$.

(iv) This follows from \refpp{r5}(iii).
\eop

\begin{rem}\labttr{r=(h-1)E8}  
Let $h$ be a fixed point free order 5 element in $O(E_8)$. Then $h$ is unique
up to conjugacy in $O(E_8)$ and $(h-1)E_8$ has the discriminant group $5^4$ and
has maximal Witt index.  In addition, $(h-1)E_8$ is rootless \refpp{r4.7}.
Thus, $(h-1)E_8\cong R$.
\end{rem}

\begin{lem}\labtt{r13} We use the notation of \refpp{r5}(iii). Let $Y=E_8$ and let $X\le Y$  be a rootless sublattice so that $Y/X \cong 5^2$.

(i)
$O(X)\cap O(Y)$ has order 2400 and is isomorphic to
$5{:}GL(2,5)\cong (5\times SL(2,5)){:}4$, order $2^53{\cdot}5^2$.

(ii)
$O(X)\cap O(Y)$  acts faithfully on $\dg X$ and on $\dg Y$ and has index $12$ in $O(X)$.

(iii)
$O(X)$ acts faithfully on $\dg X$ and induces $O(\dg X)\cong O^+(4,5)\cong (SL(2,5)\circ SL(2,5)).2^2$ on it.  In particular, $|O(X)|=2^73^25^2$.
\end{lem}
\pf  (i)  This follows from \refpp{r5}(iii).  Since $O(X)\cap O(Y)$ induces $GL(2,5)$ on $Y/X$, the kernel of the action on $\dg X$ is contained in the normal subgroup of order 5.  Such an element of order 5 in fact acts on $\dg X$ and on $\dg Y$ with a pair of Jordan blocks of degree 2 (this follows from  \refpp{r5}(iv)).
Compare with \refpp{fourdimog }.  Thus, (ii) follows.

(iii): By transitivity in \refpp{r8}(iii), we get $|O(X):O(X)\cap O(Y)|=12$, whence $|O(X)|=|O^+(4,5)|$.  By faithful action (i), we get $O(X) \cong O^+(4,5)$.
\eop

\medskip

Now Theorem \ref{mainth1} follows from Lemma \ref{r8} (ii) and Lemma
\ref{r13} (iii).

\begin{lem}\labtt{r15}  (i) There exists a sublattice in $R$ isometric to $A_4(1)\perp A_4(1)$.
We obtain an overlattice $R$ from a  sublattice of shape $A_4(1)\perp A_4(1)$ in $R$ are obtained by gluing a singular vector in $\dg {A_4(1)}$ to a singular vector in $\dg {A_4(1)}$.
Therefore, $R$ is isometric to a sublattice of index 5 in $A_4\perp A_4$.

(ii)
The cosets of $R$ in $\dual R$ have the following minimal norms

zero coset: 0

singular coset: 2

nonsingular coset (the numerator is a square mod 5):  $\frac {4}5, \frac {6}5$;

nonsingular coset (the numerator is a nonsquare mod
5):  $\frac {8}5, \frac {12}5$;
\end{lem}
\pf
(i)
For the first statement, see \refpp{r5}(iv).
A look at minimum norms in duals \refpp{r6} shows that  only a gluing of singular cosets allows the overlattice $R$ to be rootless.  Now use $O(A_4(1))\cong 2 \times PGL(2,5)$ \cite{glee8}.

(ii) By \refpp{orderog} and Witt's theorem applied to the action of $O(R)$ on
$\dg R$, the norms in $\FF_5$ indicate the orbits.
The norms may be read off from \refpp{r3} because of the gluing described in (i).
\eop

\section{Uniqueness of  $Q$}

\begin{nota}\labttr{q1plus}
We let $Q$ be a rank 16 rootless integral lattice with
$\dg Q\cong 5^4$. By Lemma \ref{maximalWitt}, $\dg Q$ is a quadratic space with
maximal Witt index. We let $U$ be a unimodular rank 16 even integral lattice
which contains $Q$. The Witt classification \cite{witt} allows two
possibilities, $U\cong HS_{16}$ or $U\cong E_8\perp E_8$.
\end{nota}

Next we shall study the isometry type of $Q$ and establish Main Theorem
\ref{mainth2}.

\begin{lem}\labtt{nonembq}  A lattice $Q$ as in
\refpp{q1plus} does not embed in $HS_{16}$.
\end{lem}
\pf Suppose that there is such an embedding of $Q$ into $S\cong HS_{16}$.
Then $S/Q \cong 5^2$.  If $T$ is the root sublattice of $S$, then $Q\cap T$ is a
rootless sublattice of $T$.  This corresponds to an elementary abelian group of
order 25 in the Lie group of type $D_{16}$ whose centralizer is just a torus.   By
\refpp{qinhs16}, no such group exists. \eop




\begin{thm}\labtt{qine8e8.3}
We let $Q$ be as in  \refpp{q1plus}.  Then $Q$ is
unique up to isometry.  In particular, $Q\cong DIH_{10}(16)$ \cite{glee8}.
\end{thm}
\pf \refpp{nonembq}, \refpp{qine8e8.2}.
\eop




\begin{lem}\labtt{q3}
Let $U_1, U_2$ be the orthogonal direct summands of $U$ which are isometric to $E_8$.

(i)
Let $W$ be any sublattice of $U$ which is rootless and has determinant $5^4$.
Then $W\cap U_i$ is isometric to $R$, for $i=1,2$.

(ii) The natural map  of $W$ in $\dg {(W\cap U_i)}$ is onto, for $i=1,2$.
\end{lem}
\pf (i) Since $W$ is rootless, so are the $W\cap U_i$ and so each $|U_i : U_i \cap
W|$ is divisible by $5^2$ \refpp{r4}.   Since $det(W)=5^4$, $5^2=|U:W|\ge |U_i
: U_i \cap W|$, for $i=1,2$, whence equality.  Now use \refpp{r8}.

(ii) This follows from (i) and determinant considerations.
\eop

\begin{rem}\labttr{q=m+n}
As in Section 2 of \cite{gl3cpath}, we can obtain an explicit embedding of $Q\cong DIH_{10}(16)$ in $E_8\perp E_8$ as follows.

Let $h\in O(E_8)$ be a fixed point free element of order $5$. Denote
\[
\begin{split}
M&:=\{(\a,\a)\in E_8\perp E_8 \mid \a\in E_8\}, \\
N&:=\{(h \a,\a)\in E_8\perp E_8 \mid \a\in E_8\}.
\end{split}
\]
By direct calculation, it is easy to show that $t_Mt_N=h^{-1}\oplus h$ has
order $5$. Notice that $t_M(\a,\b)=-(\b,\a)$ and $t_N(\a, \b)=-(h\b, h^{-1}\a)$
for any $(\a, \b)\in E_8\perp E_8$. Moreover, $M+N$ is rootless since
$(h-1)E_8$ is \refpp{r4.7}. Thus, $M+N\cong DIH_{10}(16)\cong Q$ by the
classifcation of rootless $EE_8$ pairs \cite{glee8}.
\end{rem}



\begin{lem}\labtt{q7}
(i) Let $J$ be an overlattice of $R_1 \perp R_2$, where $R_1\cong
R_2\cong R$, such that $J$ is obtained by gluing 2-spaces
$W_1$ from $\dg {R_1}$ and $W_2$ from $\dg {R_2}$. If $J$ is integral, then $W_1$ and $W_2$ are isometric.

(ii)
An integral
overlattice $J$ of $R_1 \perp R_2$ corresponding to a gluing of 2-spaces from orbits of type  $Orbit(a,\cdots )$ (see \refpp{q6}) in the respective quadratic spaces
$\dg {R_1}$ and $\dg {R_2}$ is denoted
$L(a,\cdots )$.  The minimum norms are given below.
The associated discriminant groups have maximal Witt index.

$L(2)$: 4,

$L(1,s)$: 2,

$L(1,n)$: 4,

$L(0,1)$:  2,

$L(0,0)$: 2.
\end{lem}

\pf (i) The projections $p_i$ to the spaces $\QQ\otimes R_i$ have the property that for all $v\in J$, $p_i(v)$ has norm in $\frac 25 \ZZ$ for $i=1,2$ and that
the sum of these two norms is in $2\ZZ$.  This means that the associated linear isomorphism $W_1\rightarrow W_2$ is the negative of an isometry.  Since $-1$ is a square modulo 5, the quadratic spaces $W_1, W_2$ are isometric.


(ii) Use \refpp{q6}, \refpp{r15} and the fact that $-1$ is a square modulo 5.
\eop

\begin{coro}\labtt{q7.5}
We use the notation of \refpp{q7}.
Suppose $J$ is isometric to $Q$ and is in $Orbit(2)$.
Let $p_i$ be the orthogonal projection to $\QQ \otimes R_i$, $i=1, 2$.
Then

(i) $Stab_{O(R_i)}(p_i(J))$ has shape $5{:}GL(2,5)\cong (5\times SL(2,5)){:}4$, for $i=1, 2$.

(ii) The kernel of the action of
$Stab_{O(R_i)}(J)$  on $p_i(J)/R_i$ is the normal subgroup of order 5,  for $i=1, 2$.

(iii)
 For $i=1,2$, $Stab_{O(J)}(R_1)=Stab_{O(J)}(R_2)$ and this group may be interpreted as a fiber product for the maps
$Stab_{O(R_i)}(J) \rightarrow GL(2,5)$, $i=1,2$, i.e.,
the set of pairs
$(g_1, g_2)\in Stab_{O(R_1)}(J) \times Stab_{O(R_2)}(J)$ such that
the images in $GL(2,5)$
of $g_1$ and $g_2$ are equal.


(iv) $O(J)\cap O(R_1\perp R_2)$
has shape $(5\times 5)(GL(2,5) \times 2) \cong (Dih_{10} \times 5){:}GL(2,5)$.
The stabilizer in this group of $R_1$ or $R_2$ has index 2.
The
kernel of the action on $\dg J$ is a group of order 10, generated by
$EE_8$-involutions.  Let $t$ be any $EE_8$-involution in $O(J)\cap O(R_1\perp R_2)$.
We have $O(J)\cap O(R_1\perp R_2) = C_{O(J)\cap O(R_1\perp R_2)}(t)F$, where $F$ is the normal subgroup of order 5 inverted by $t$.

\end{coro}
\pf  (i)  This follows from the structure of the stabilizer in $O^+(4,q)$ of a
maximal totally singular subspace.   This stabilizer has shape
$TH$, where $T$ is normal, $T=U \times V$, where $U$ is isomorphic to $\FF_q$, $U\cong SL(2,q)$ and where $H\cong \ZZ_{q-1}$, $UH\cong AGL(1,q)$ and $VH\cong GL(2,q)$.

(ii) \refpp{r5}(iii).

(iii) This follows from (i) and (ii).

(iv) Use (i,ii,iii) and the facts that $O(J)\cap O(R_1\perp R_2)$
has an $EE_8$-involution which interchanges $R_1$ and $R_2$ and contains
$Stab_{O(J)}(R_1)=Stab_{O(J)}(R_2)$ with index 2.   Note \refpp{trivialaction}.
\eop

\begin{coro}\labtt{q7.6}
We use the notation of \refpp{q7}.
Suppose $J$ is isometric to $Q$ and is in $Orbit(1,n)$. Then

(i) $Stab_{O(J)}(R_i)$ has shape $5^2(4\times 2 \times 2)$.

(ii) The kernel of the action of
$Stab_{O(R_i)}(J)$  on $J/(R_1 \perp R_2)$ is a normal subgroup of order 5.

(iii)  $O(J)\cap O(R_1\perp R_2)$ is isomorphic to
$Dih_{10} \times 5.(4\times 2 \times 2)$. The
kernel of the action on $\dg J$ is a group of order 10, generated by
$EE_8$-involutions.
\end{coro}
\pf
Use \refpp{q6.5} and follow the arguments of \refpp{q7.5}.
\eop

\begin{lem} \labtt{q7.8} $|O(Q)|=10|O^+(4,5)|=2^{8}3^25^2$ and $O(Q)$ induces $O^+(4,5)$ on $\dg Q$.
The kernel of the action is $K$, the dihedral group generated by a pair of $EE_8$- involutions whose negative eigenspaces generate $Q$.
\end{lem}
\pf The order of $O(Q)$ is determined since $O(Q)$ acts transitively on the  set of
overlattices isometric to $E_8\perp E_8$  and the stabilizer of one of them has
the form given in \refpp{q7}(v). \eop

\begin{coro}\labtt{q8.5}
There exists an embedding
$O(Q)\rightarrow Frob(20) \times O^+(4,5)$ such that the image has index 2 and contains neither direct factor.
\end{coro}
\pf
A consequence of \refpp{autdihee8} is that $O(Q)$ maps onto $Aut(Dih_{10})\cong Frob(20)$ by its doubly transitive action on the involutions of $\la t, u \ra$.
Therefore, the normal subgroup $H$ of $O(Q)$ which centralizes $\la t , u \ra$ has index 20.  It follows that $H$ acts faithfully on $\dg Q$ and induces on $\dg Q$ a subgroup of index 2 in $O^+(4,5)$.  We have
$H\cap \la t , u \ra = 1$.
\eop

\begin{lem}\labtt{q4}
Since $Q\cong DIH_{10}(16)$, $Q$ is spanned by a pair of $EE_8$-lattices.  The
SSD-involutions corresponding to any $EE_8$-sublattice stabilize any $E_8^2$
overlattice of $Q$ and interchange its two indecomposable  summands.
\end{lem}

\pf Use \refpp{rssdonsum}.
 \eop


\begin{lem}\labtt{QinQ*}
The cosets of $Q$ in $Q^*$ have the following minimal norms:

zero coset: 0

singular coset: 2

nonsingular coset (the numerator is a square mod 5) :  $\frac {14}5, \frac {16}5$;

nonsingular coset (the numerator is a nonsquare mod
5) :  $\frac {8}5, \frac {12}5$;
\end{lem}
\pf Let
$U$ be an overlattice of $Q$, where
$U=U_1\perp U_2$, $U_i\cong E_8, i=1,2$.

Denote $R_i=Q\cap U_i$, $i=1,2$. Then $R_1 \cong R_2\cong R$.  We may
also assume
\[
Q=span \big\{\{(\a, \a)\mid \a\in E_8\}\cup R_1\cup R_2\big\}.
\]
Then
\[
Q^* =span\big\{  \{(\gamma, -\gamma) \mid \gamma\in R^*\}\cup U_1\cup U_2\big\}.
\]
The minimum  norms of the cosets may be read off from \refpp{r15}. \qed


\medskip

We summarize our results on rank 16 lattices.

\begin{nota}\labttr{apm}
We define the following pairs of lattices in Euclidean space $\RR^{16}$:  

$\mathcal A := \{ (X,Y) \mid X\le Y, X\cong R_1 \perp R_2, Y\cong Q\}$,

$\mathcal A^+ :=\{ (X,Y)\in {\mathcal A} \mid Y \in Orbit(2) \text{ with respect to
}X; \text{ see }\refpp{q7} \}$ and

$\mathcal A^- :=\{ (X,Y)\in {\mathcal A} \mid Y \in Orbit(1,n) \text{ with respect
to }X; \text{ see }\refpp{q7}\}$.

$\mathcal B : = \{ (Y, Z) \mid  Y\le Z,
Y \cong Q, Z\cong E_8\perp E_8 \}$.
\end{nota}

\begin{thm}\labtt{q9}
(i)  The group $O(16,\RR )$ acts transitively on each of the sets $\mathcal A^+,
\mathcal A^-, \mathcal B$.

\noindent

$(ii.a+)$ If  $(X,Y)\in \mathcal A^+$, we have
\begin{align*}
&O(X) \cong
O^+(4,5)\wr 2, &  &|O(X)|=2^{15} 3^45^4, \cr
&O(Y) \cong Dih_{10}{\cdot}O^+(4,5), & &
|O(Y)|=2^8 3^2 5^3, \cr
&O(X)\cap O(Y) \cong (Dih_{10} \times 5){:}GL(2,5), & &
|O(X)\cap O(Y) |=  2^6 3{\cdot}  5^2,\cr
&|O(X):O(X)\cap O(Y)|= 2^{9}3^35^2,  & &
|O(Y):O(X)\cap O(Y)|=  2^2 3{\cdot} 5 .
\end{align*}

$(ii.a-)$ If  $(X,Y)\in \mathcal A^-$, we have
\begin{align*}
& O(X) \cong
O^+(4,5)\wr 2, & &|O(X)|=2^{15} 3^45^4,\cr
&O(Y) \cong Dih_{10}{\cdot}O^+(4,5) ,& &
|O(Y)|=2^8 3^2 5^3, \cr
&O(X)\cap O(Y) \cong 5^2{:}(4 \times 2 \times 2 \times 2), &&
|O(X)\cap O(Y) |= 2^5 5^2 ,\cr
&|O(X):O(X)\cap O(Y)|= 2^{10}3^45^2,  &&
|O(Y):O(X)\cap O(Y)|=  2^3 3^2 5.
\end{align*}


(ii.b) If  $(Y,Z)\in \mathcal B$, we have
\begin{align*}
& O(Y) \cong
Dih_{10}{\cdot}O^+(4,5), && |O(Y)|=2^8 3^2 5^3, \cr
&O(Z) \cong O(E_8)\wr 2, & &
|O(Z)|=2^{29} 3^{10} 5^4 7^2, \cr
&O(Y)\cap O(Z) \cong Dih_{10}.(5\times GL(2,5)),
&& |O(Y)\cap O(Z) |= 2^6 3{\cdot}5^3,  \cr
&|O(Y):O(Y)\cap O(Z)| =  2^2 3,   & &
|O(Z):O(Y)\cap O(Z)|=  2^{23}3^9 5{\cdot }7^2.
\end{align*}


(iii)
 Given $Y \cong Q$, there is a bijection between
the set of $X$ such that $(X,Y)\in \mathcal A^+$ and the set of $Z$ such that
$(Y,Z)\in \mathcal B$.  The bijection may be described as follows.

If we write $X=X_1\perp X_2$ ($X_1\cong X_2 \cong R$), the associated $Z$ is $proj_{X_1}(Y) \perp proj_{X_2}(Y)$.

If we write $Z=Z_1\perp Z_2$ ($Z_1\cong  Z_2 \cong E_8$), the associated $X$ is $(Z_1\cap Y)\perp (Z_2\cap Y)$.
\end{thm}
\pf
For orders of $O(R), O(Q)$, see \refpp{r13}, \refpp{q7.8}.  The order of $O(E_8)$ is $2^{14}3^55^27$.

For $\mathcal A^+$, see
 \refpp{q7.5}(iii).

For $\mathcal A^-$, see
 \refpp{q7.6}(iii).

For $\mathcal B$, see \refpp{q7.8} and note that an even unimodular overlattice of $Q$ must correspond to a maximal totally singular subspace of $\dg Q$.
\eop

\section{Embeddings  of $Q$ and $R$ in Niemeier lattices}

\begin{nota}\labttr{el1}
We suppose that $\L$ is a Leech lattice containing a sublattice $Q'$  isometric to
$Q$. We define $R':=ann_{\L}(Q')$.  Then $R'$ is an even rank 8 lattice and $\dg
{R'} \cong \dg {Q'} \cong 5^4$ \cite{gal}.  By \refpp{r8}, $R'\cong R$.
\end{nota}





\begin{prop}\labtt{el3}
Let $N$ be a Niemeier overlattice of $Q\perp R$. Let $X:= N\cap R^*$ and
$Y:=N\cap Q^*$.  Define the integer $e$ by $5^e=|X:R|=|Y:Q|$.  Then $e\in \{0,1,2\}$.  There exists a Niemeier overlattice of $Q\perp R$ which realizes each of these values of $e$.

(i) If $e=0$, then $N\cong \L$. Moreover,
$O(Q\perp R)$ acts transitively on the set of Niemeier overlattices which contain
each of $Q$ and $R$ as a direct summand. Therefore, the embedding of $Q\perp
R$ into $\L$ is essentially unique;

(ii) If $e=1$, then $N\cong \mathcal{N}(A_4^6)$ and
$O(Q\perp R)$ acts transitively on the set of Niemeier overlattices such that
$X/R\cong Y/Q\cong 5$; these Niemeier lattices are isometric to $\mathcal N(A_4^6)$.

(iii) If $X/R\cong Y/Q\cong 5^2$, then $N\cong E_8^3$ and   $O(Q\perp R)$  acts
transitively on the set of Niemeier overlattices such that $X/R\cong Y/Q\cong
5^2$;  these Niemeier lattices are isometric to $E_8^3$. 
\end{prop}
\pf Since an integral overlattice corresponds to a totally singular subspace of
the discriminant group, $e\le 2$. In all three cases, transitivity follows from
\refpp{el1.1}, \refpp{el2} and the fact that the projection of $N$ to $\dual R$
lies in $\dual X$.
Finally, we must settle the isometry types in those
respective orbits.

\medskip




Since the minimal norms for $Q^*$ and $R^*$  are $8/5$ and $4/5$, respectively
(see \refpp{QinQ*} and \refpp{r15}) and $8/5+4/5 =12/5 > 2$, the roots of $N$
must be  in $Y\perp X$.   For case (i), this means $N$ is rootless, so is
isometric to the Leech lattice.   For (ii), note that since $X$ is a direct
summand, the root sublattice of $N$ has an indecomposable orthogonal component
isometric to $A_4$. Since the Coxeter number of all components of the root
system is the same \cite{venkov}, the root system of $N$ is $A_4^6$ and $N\cong
\mathcal{N}(A_4^6)$.

\medskip

(iii) If $X/R\cong Y/Q \cong 5^2$, then $X$  and $Y$ are both even unimodular.
Thus,  $X\cong E_8$ and $Y\cong E_8^2$ or $HS_{16}$. However,  $Q$ cannot
be embedded into $HS_{16}$  by \refpp{nonembq} and hence $N\ncong
\mathcal {N}(D_{16} E_8)$.  Therefore,  $N\cong E_8^3$.
\eop

\medskip

 Main Theorems \ref{mainth3} and \ref{mainth4} now follow from
\refpp{q8.5} and \refpp{el3}.

\medskip

Finally, we shall prove Corollary \ref{maincoro1} and Main Theorem
\ref{mainth5}.

\medskip

\noindent\textbf{Proof of Corollary \ref{maincoro1}.}  
The centralizer in $O(\L)$ of $D$  is isomorphic to a group $H$ of index 2 in
$O^+(4,5)$ \refpp{q8.5}, \refpp{el3}. 
Let $Z:=Z(O^+(4,5))$. We have that $H/Z$ is isomorphic to one of
$Alt_5 \wr 2$ or the even subgroup of $Sym_5 \times Sym_5$
\refpp{fourdimog2}. Both of these groups contain subgroups isomorphic to
upwards extensions of $PSL(2,5)$ by nonidentity cyclic 2-groups but only the
latter contains such a subgroup which normalizes a group of order 5. By
\refpp{q7.5}, $Alt_5 \wr 2$ can not occur as the central quotient of
$C_{O(\L)}(D)$. \eop

\medskip


\medskip

\noindent\textbf{Proof of Main Theorem \ref{mainth5}.}

We use \refpp{el3} and \refpp{staboverlattice}(iii), applied to
$X$ and $Y$ in the notation of \refpp{el3} to obtain a group $S$.  The group $H$
is $S\cap O(Q\cap R)$.

For all $e=0, 1, 2$,
$S\cap O(Q\cap R)$ covers $S/(S\cap O(X) \times (S\cap O(Y))$ since
$O(Q), O(R)$ induces the full orthogonal group on $\dg Q, \dg R$, respectively; see \refpp{el2}.
Therefore, $H$ is an extension of normal subgroup
$\{ g \in O(R)\cap O(X) \mid (g-1)N\le X \} \times
\{ g \in O(Q)\cap O(Y) \mid (g-1)N\le Y \}$ by
quotient isomorphic to $O(\dg X)\cong O(\dg Y)$ (we have $\dg X\cong \dg Y\cong \FF_2^{4-2e}$, a quadratic space of maximal Witt index).
We give descriptions of the normal subgroup
$\{ g \in O(R)\cap O(X) \mid (g-1)N\le X \} \times
\{ g \in O(Q)\cap O(Y) \mid (g-1)N\le Y \}$ in the respective cases below.


For $e=0$, $\{ g \in O(R)\cap O(X) \mid (g-1)N\le X \} \cong 1$ and
$\{ g \in O(Q)\cap O(Y) \mid (g-1)N\le Y \} \cong Dih_{10}$.

For $e=1$,
$\{ g \in O(R)\cap O(X) \mid (g-1)N\le X \} \cong 5$ and
$\{ g \in O(Q)\cap O(Y) \mid (g-1)N\le Y \} \cong Dih_{10} \times 5$.

For $e=2$,
$\{ g \in O(R)\cap O(X) \mid (g-1)N\le X \} = O(R)\cap O(X) \cong 5{:}GL(2,5)$ and
$\{ g \in O(Q)\cap O(Y) \mid (g-1)N\le Y \} = O(Q) \cap O(Y) \cong (Dih_{10} \times 5){:}GL(2,5)$.

\section{A suggestion of triality}

The moonshine VOA \cite{FLM} was built from a VOA
based on the Leech lattice and some of its modules.  A so-called extra automorphism (sometimes called ``triality'') was constructed which, with a natural group acting on the VOA, generated a copy of the Monster acting as VOA automorphisms.

For the $3C$-case \cite{gl3cpath}, a  Weyl group is associated to one type of overlattice and triality is used to create a twist of this Weyl group and a corresponding twist of glue map, which is used to define a
different type of overlattice.  (The Weyl group of $E_8$ has a unique nonsolvable compostion factor isomorphic to $D_4(2)$, which has a group of graph automorphisms representing triality.)  The ``loss'' of half the Weyl group was due to passing to a half-spin module which affords a projective representation of the special orthogonal group but not the whole orthogonal group.

\begin{rem}\labttr{updownpsl25psl25}
For our $5A$-situation, it is not obvious how to imitate the above program.
Nevertheless, one can point to weak analogues of our $3C$ story.
For example,   the cases $e=0$ and $e=2$ in \refpp{mainth5} give a central extension of $SL(2,5)\circ SL(2,5)$ to $SL(2,5)\times SL(2,5)$.  Such an extension
can be created by lifting from the special orthogonal group
to the spin group (of type $D_4$, which has triality).

If we consider the isometry groups of the three kinds of Niemeier overlattices, we
see a variety of upward and downward extensions of $PSL(2,5) \times PSL(2,5)$ within the centralizer of the dihedral group of order 10.

Let $N$ be a Niemeier overlattice and $D$ a dihedral group of order 10 generated by a pair of $EE_8$-involutions.

If $N\cong E_8^3$, $C_{O(N)}(D)\cong O(E_8) \times Dih_{10}\times 5 \times 2$;
this group contains both $Sym_5\wr 2$ and $O^+(4,5)$.

If $N\cong \L$, $C_{O(N)}(D)\cong (SL(2,5)\circ SL(2,5)){:}2$;  this group does not contain a subgroup isomorphic to $Alt_5 \times Alt_5$.

If $N\cong {\mathcal N}(A_4^6)$,
$C_{O(N)}(D)$
contains $Alt_5 \times Alt_5$ but not $SL(2,5)\circ SL(2,5)$
\refpp{centdniema4to6}.

\end{rem}

\section{Pieces of Eight, \`a la  $DIH_{10}(16)$}

The basic POE program \cite{poe} to study the Leech lattice, Mathieu and Conway groups, was centered on a study of $EE_8^3$-lattices.  This program can be carried out using a new viewpoint.  We start with a triple of pairwise orthogonal lattices $R_1, R_2, R_3$ which are isometric to $R$, classify overlattices of $R_1 \perp R_2$  which are isometric to $Q$ and correspond to $Orbit(2)$ in the sense of \refpp{q6}, then finally glue such a $Q$ with $R_3$ to get a Leech lattice.
This program offers new proofs of 5-local information about $O(\L )$.

\bigskip

\newpage

\centerline{\Large \bf  Appendices}

\appendix

\section{Some general results about lattices}

\begin{lem}\labtt{r0.5}  Let $L$ be a lattice of rank $p-1$ and $h$ an isometry of order $p$.  Suppose that $M, N$ are $h$-submodules of $\QQ L$ and that
$N$ has index a power of $p$ in $M$.  Then $M/N$ is a uniserial module.
\end{lem}
\pf
The determinant of $h-1$ is $p$.  Therefore the series
$(h-1)^i M$ for $i=0,1,\dots r$ forms a chain from $M$ down to $N$ such that each term has index $p$ in the previous.  Since $h$ has order $p$ and $M/N$ is a finite $p$-group, there are no other $h$-submodules of $M$ which contain $N$.
\eop

\begin{lem}\labtt{r0}
Suppose that $p$ is an odd prime
and that $L$ is a rank $p-1$ integral lattice with an automorphism $h$ of order $p$.
Suppose that $\dg L$ is elementary abelian $p$-group of rank $r\in \{1, 2, \dots, p-1\}$.
Let $i, j \in \{0, 1, \dots , r\}$.
Then $((h-1)^i x, (h-1)^j y )=0$ for all $x, y \in L$ if and only if $i+j\ge r$.
\end{lem}
\pf
Consider the descending chain $(h-1)^i \dg  L$ for $ i\ge 0$ and its
ascending chain of annihilators \refpp{r0.5}.
\eop

\begin{lem}\labtt{r1}
Suppose that $p$ is an odd prime
and that $L$ is a rank $p-1$ integral lattice with an automorphism $h$ of order $p$.  Suppose that $(L,L)\not \le p\ZZ$ and that $\dg L$ is elementary abelian $p$-group of rank $r\in \{0,1,\dots, p-1\}$.
Then the nonsingular quadratic space $L/p\dual L$ is isometric to
$(h-1)L/p\dual {((h-1)L)} \perp U$, where $U$ has Gram matrix $\frac 1p \begin{pmatrix}0&1\cr 1&0 \end{pmatrix}(mod\,  \ZZ )$.
\end{lem}
\pf
The commutator space $(h-1)^{r-1}\dg L$ is 1-dimensional and consists of singular vectors.  Let $W$ complement $(h-1)^{r-1}\dg L$ in $ann((h-1)^{r-1}\dg L)=(h-1)\dg L$.
Then $W$ is nonsingular and $ann(W)$ is a 2-dimensional
space containing a singular vector, so is split.
By \refpp{r0.5} and determinant considerations,  $\dual {((h-1)L)}$ contains $\dual L$ with index $p$ and
$p\dual {((h-1)L)}/pL$ is the 1-space $(h-1)^{r-1}\dg L$.
\eop

\begin{coro}\labtt{r2}
The discriminant group of $A_4(1)$ is isometric to the orthogonal direct sum of the  discriminant group of $A_4$ and the $\FF_5$-space with
 Gram matrix  $\frac 15 \begin{pmatrix}0&1\cr 1&0 \end{pmatrix}(mod\,  \ZZ )$.    Also, the minimum norm in $\dual {A_4}$ is $\frac 45$.
\end{coro}
\pf \refpp{r1}.  \eop

\begin{lem}\labtt{r3}
Let $h$ be an isometry of order 5 on $A_4$.  Define $A(i):=(h-1)^iA_4$, for $i \in \ZZ$.  The lattices in the chain below have the indicated minimum norms (note that $A_4(1)\cong (h-1)A_4$ and that if $X$ is an $h$-invariant lattice in $\QQ A_4$, then $(h-1)^2X\cong \sqrt 5\, X$ and $(h-1)^4X=5X$).

$\begin{matrix}
\cdots & A(4) & A(3)&A(2)&A(1)&A(0)&A(-1)     &A(-2)   &A(-3)         &A(-4)\cdots            \cr
\cdots & 50   &20    &10   &4     &2     &\frac 45 &\frac 25 &\frac 4{25} & \frac 2{25} \cdots \cr
\end{matrix}  $
\end{lem}
\pf
Note that $h-1$ has determinant 5 and so
$|X:(h-1)X|=5$.
Also, $(h-1)^i$ has determinant $5^i$ and $(h-1)^4$ induces the 0-map on
$X/5X$, so $(h-1)^4X=5X$.

The adjoint of $h$ is $h^{-1}$ and so
$((h-1)x,(h-1)y)=(x,(h^{-1}-1)(h-1)y)=(x,(h-1)^2(-h^{-1}y)$.
It follows that $A(2)$ satisfies
$(A(2),A(2))\le 5\ZZ$.  Therefore
$\frac 1{\sqrt 5}A(2)$ is an even integral lattice of determinant 5.  By \cite{glee8,gal}, this lattice is isometric to $A_4$.

Finally, recall that $A_4(1)$ is isometric to $(h-1)A_4$ \cite{glee8,gal}.
Since $A(3)=(h-1)A(2)$, the previous paragraph and the definition of $A_4(1)$ and $(h-1)A_4$ implies that $A(3)\cong \sqrt 5 A_4(1)$.

The first statement follows.  The table of minimum norms is a consequence of the known minimum norms in $A_4$ and $A_4(1)$, respectively.
\eop

\begin{lem}\labtt{r4.5}
Suppose that $L$ is a lattice invariant under an isometry $g$ of prime order $p$ and that the minimal polynomial of $g$ is the cyclotomic polynomial
$(x^p-1)/(x-1)$.

(i) Let $n=rank(L)$.  There exists an integer $m$ so that $n=(p-1)m$.

(ii) If $K$ is  a rank $n$ $g$-invariant lattice in $\QQ \otimes L$, then
$(g-1)^{p-1}K=pK$ and
the action of $g$ on $K/pK$ has Jordan canonical form the sum of $m$ blocks of degree $p-1$.
\end{lem}
\pf
By \cite{rim}, we may write
$L$ as a direct sum of indecomposable modules $L_i$, $i=1,\dots, q$,  for $g$.  Each $L_i$ is isomorphic to an ideal in $\ZZ [\la g \ra]/(1+g+\cdots + g^{p-1})$ \cite{rim} and so has rank $p-1$ as an abelian group and on it the minimal polynomial of $g$ is $(x^p-1)/(x-1)$.  Therefore $q=m$ and the elementary divisors for $g-1$ on $L$ are $(1, \dots , 1, p, \dots p)$, where $p$ has multiplicity $m$.  It follows that
$K/(g-1)K \cong p^m$ (elementary abelian).  Since
$(g-1)^{p-1} \equiv 1+g+\cdots +g^{p-1}\equiv 0 (mod\, p)$,
$(g-1)^{p-1}(K)\le pK$.  Since $|K : (g-1)^{p-1}K|=p^{m(p-1)}=p^n=|K:pK|$, we have
$(g-1)^{p-1}(K)= pK$.
Each Jordan block for the action of $p-1$ on $K/pK$ has degree at most $p-1$ and if there exists a block of size less than $p-1$, there are at least $m+1$ nontrivial blocks, whence $K/(p-1)K$ has dimension at least $m+1$, a contradiction.
\eop

\begin{lem} \labtt{r14}
(i)
The minimum norm for glue vectors of $A_4$ in $\dual {A_4}$ is $\frac 45$;
the nontrivial cosets have minimum norms
$\frac 45$ and $\frac 65$.

(ii)
Norms of vectors in the lattice $\dual {A_4(1)}$ have the form $\frac k5$ where $k$ is an even integer.
The minimum norm for nonzero glue vectors of $A_4(1)$ in $\dual {A_4(1)}$ are  $\frac 25, \frac 45, \frac 65, \frac 85, 2$ \cite{gal}.
\end{lem}
\pf
(i)
A vector in $\dual A_4 \setminus A_4$ has norm of the form $\frac k5$, where $k$ is an even integer.
The vector  $\pm \frac 15 (1,1,1,1,-4)$ in the standard model of $A_4$  has norm
$\frac 45$ and it is easy to prove that there is no smaller norm and that
a minimum norm vector in $\dual A_4 \setminus A_4$ is just this vector, up to negation and coordinate permutation.  For details, see \cite{gal}.  For if $u\in \dual {A_4}$ has the minimum, $u$ has form $\pm \frac 15 (1,1,1,1,-4)$ in the standard model of $A_4$.
The coset
$2u+A_4$ contains the vector $\pm \frac 15 (2,2,2,-3,-3)$, which has norm $\frac 65$ and no smaller norm is possible in this coset.

(ii)
The first statement follows from the fact that $\dg {A_4(1)}$ is an elementary abelian group of exponent 5 and that $A_4(1)$ is an even lattice.

For a singular coset, the minimum is norm 2 since $A_4$ is properly between $A_4(1)$ and $\dual {A_4(1)}$.
The vectors mentioned in (i) and in $A_4$ are in the dual of $A_4(1)$.
Therefore, minimum norms $\frac 45, \frac 65, 2$ occur.  It remains to deal with the cosets
of $A_4(1)$ in $\dual {A_4(1)}$ where  the norm of a vector has the form  $\frac k5$ where $k \in \pm 2+10\ZZ$.

In the notation of \refpp{r3}, $\dual {A_4(1)}=A(-2)$ and we get minimum norm
in $A(-2)$ is $\frac 25$.  If $v$ is such a minimal vector, we have $(2v,2v)=\frac
85$ and no vector in $2v+A_4(1)$ has lesser norm since $|\frac 85|<2$. \eop


\begin{lem}\labtt{r9.5}
For an integer $k$ and a finite abelian group $A$, denote by $A_{k'}$ the subgroup of all elements of $A$ whose order is a number relatively prime to $k$.

Let $L$ be an integral lattice and $M$ a sublattice of finite index $n$.
The natural map
$\dual L\rightarrow \dual M \rightarrow \dg M \rightarrow {\dg M}_{n'}$ induces an isometry
${\dg L}_{n'}\cong {\dg M}_{n'}$.
\end{lem}
\pf
We have the containments
$M\le L \le \dual L \le \dual M$.   Since $n=|L:M|=|\dual M : L|$, the natural map $\varphi$
of $\dual M$ onto ${\dg M}_{n'}$ remains onto when restricted to $\dual L$.
If $\varphi |_{\dual L}$ has kernel $K$, then $K/M$ has the property that every element of it is annihilated by a power of $n$.  Since $\dual L/K$ has order prime to $n$, we have a splitting $\dual L/M=K/M \oplus J/M$, for a sublattice
$J$, where $M\le J \le \dual L$ and $(|J/M|,n)=1$.   Clearly, $J$ is uniquely determined and $J$ maps isomorphically onto ${\dg M}_{n'}$.  Since $K\cap J=M$, $J/M\cong {\dg M}_{n'}$, as required.    Since $M \le J \le \dual M$, the isomorphism of groups $J/M\cong {\dg M}_{n'}$ is an isometry.
\eop


\begin{lem}\labtt{trivialaction}
Let $L$ be an integral lattice and $M$ an SSD sublattice.  The involution associated to $M$ acts trivially on $\dg L$.
\end{lem}
\pf
Let $P^{\vep}$ be the orthogonal projection to the
$\vep$-eigenspace of $t=t_M$, $\vep = \pm $.
For $x\in L$, $x=P^+(x)+P^-(x)$ and $tx=P^+(x)-P^-(x)$.
Therefore, $x-tx=2P^-(x)$.

Now take $x\in \dual L$.
Since $L$ is integral, $P^-(x)\in \dual M \le \half M$.  Therefore,
$x-tx=2P^-(x)\in M\le L$, proving the result.  \eop

\begin{prop}\labtt{autdihee8} (i)
Suppose that we have a pair of $EE_8$-lattices $M, N$.
Let $t, u$ be the associated involutions.  Suppose that $t', u'$ are involutions of
$\la t, u \ra$ such that both $t' u'$ are associated to $EE_8$-sublattices $M', N'$ and that
$\la t, u \ra = \la t', u' \ra$.   Then there exists an isometry
$g\in O(M+N)$ so that $t^g=t'$ and $u^g=u'$.

(ii) This property (i) holds whenever $t' ,u'$ are generators and $n>1$ is an odd prime.
\end{prop}
\pf This follows from the main theorem of \cite{glee8}.
\eop

\begin{prop}\labtt{rssdonsum}
Suppose that $L=L_1\perp  L_2$ where the $L_i$ are orthogonally indecomposable lattices.  Let $M$ be a nonzero RSSD sublattice of $L$ which is a direct summand of $L$ as an abelian group and let $t=t_M$ be the associated involution.

(i) If $t$ fixes some $L_i$, then $t$ fixes both $L_i$ and
$M=M\cap L_1 \perp M \cap L_2$.

(ii) If $M$ is orthogonally indecomposable and $t$ fixes each $L_i$, then $t$ is $1$ on one of the $L_i$  and $M$ is contained in the other.

(iii) If each $L_i$ has roots, $rank(L_1)=rank(L_2)=rank(M)$ and $M$ is rootless, then $t$ does not fix any of the $L_i$.
\end{prop}
\pf
(i): This is clear since $O(L)$ permutes the set of indecomposable orthogonal summands.

(ii): This is immediate from (i).

(iii): From (ii) we get $M$ contained in $L_j$ for some $j$.  Since $M$ is a direct summand of $L$, it is a direct summand of $L_j$.   Since $M$ and $L_j$ have the same rank, $M=L_j$.  This is a contradiction since $L_j$ has roots and $M$ does not.
\eop

\begin{lem}\labtt{el1.1}
Suppose that $X, Y$ are integral lattices such that

(a) $\dg X \cong \dg Y$ has squarefree exponent; and

(b)  $O(X)$ acts on $\dg X$ as the full orthogonal group of the quadratic space $\dg X$.

Then for each divisor $d > 0$ of $det(X)$, $O(X)$ has one orbit on
$\mathcal O := \{ J \mid X\perp Y \le J \le \dual X \perp \dual Y, J \text{ is integral and }|J:X\perp Y|=d \}$, provided that this set is nonempty.
\end{lem}
\pf
A finite quadratic space is the orthogonal direct sum of its $p$-primary parts, for all primes $p$.  If $J \in \mathcal O$, $J/(X\perp Y)$ corresponds to a totally singular subspace of $\dg {X\perp Y}$.

If $W:=\dg {X\perp Y}$ and $W=\oplus_p W_p$ is the primary decomposition,
we have $O(W)\cong \prod_p O(W_p)$.  Each $W_p$ is a vector space.
Now use Witt's theorem on extensions of isometries.
\eop

\begin{lem}\labtt{el2}
Let $V$ be a quadratic space which is finite dimensional over a field and split.
 Let $W$ be a totally singular subspace.
 Then $Stab_{O(V)}(W)$ induces $O(W/Rad(W))$ on
$W/Rad(W)$.
\end{lem}
\pf
Witt's theorem  on extensions of isometries.
\eop

\begin{lem} \label{N*}
Let $M$, $N$ be sublattices of the lattice $Q$ where $N$ is a direct summand of $Q$, $Q=M+N$ and
$(det(N), det(Q))=1$.  Then the  natural map of $M$ to $\dg N$ is onto.
\end{lem}
\pf
By an elementary lemma, the natural map of $\dual Q$ to $\dual N$ is onto (for example, \cite{gre8} or  \cite{gal}.  The image under this map of $Q$ has index relatively prime to $det(N)$.  Consequently, the natural map $\phi$ of $Q$ to $\dg N$ is onto.  Since $Ker (\phi )=N\perp ann_Q(N)$ and $Q=M+N$, the natural map of $M$ to $\dg N$ is onto.
\eop

\begin{lem}\labtt{5one8}
Let $L$ be the $E_8$-lattice and $M$ and $A_4^2$-sublattice.  Let $M=M_1\perp
M_2$ be the orthogonal decomposition of $M$ into indecomposable summands.  We
let $W_1\times W_2$ be the Weyl group of $M$, where $W_i$ is the Weyl group on
$M_i$.

(i)
If $h \in W_i$ has order 5, the action of $h$ on $L/5L$ has Jordan canonical form the sum of a single degree 5 block and three degree 1 blocks.

(ii) If $h \in W_1\times W_2$ has order 5 and $h \notin W_1 \cup W_2$, then
the action of $h$ on $L/5L$ has Jordan canonical form the sum of a single degree 5 block and a single  degree 3 block.
\end{lem}
\pf
(i) Suppose $h\in W_1$.
This follows since when $v\in L \setminus M_1\perp M_2$, then $\ZZ \la h \ra v$ is a free module for
$\ZZ \la h \ra$ and it maps in $L/5L$ to a free submodule, which is an injective module.  Since $h$ fixes the rank 4 module $M_2$ pointwise and $M_2$ is a direct summand of $L$, (i) follows.

(ii)
The actions of $h$ on $L/M_2$ and on $M_2$ have minimum polynomial $x^4+x^3+x^2+x+1$.  Therefore, the action of $h$ on $L/5L$ has fixed point dimension at most 2.
As in (i), we get a 5-dimensional free module in $L/5L$ and (ii) follows.
 \eop

\begin{lem}\labtt{staboverlattice}
Suppose that $Z$ is a lattice which contains orthogonal
sublattices $X$ and $Y$ such that $X\perp Y$ has finite index in $Z$
and both $X$ and $Y$ are direct summands of $Z$.


We interpret $O(X)$ as a subgroup of $O(\QQ \otimes Z)$ by extending with trivial action on $ann(X)$.
Similarly,
we interpret $O(Y)$ as a subgroup of $O(\QQ \otimes Z)$ by extending with trivial action on $ann(Y)$.


Define $S:=Stab_{O(X)\times  O(Y)}(Z)$
and define
$T$ to be the kernel of the action of $Stab_{O(X)}(E)$ on $E$ and
$U$ to be the kernel of the action of $Stab_{O(Y)}(F)$ on $F$.
These are normal subgroups of $S$.

The overlattice $Z$ corresponds to the following data: subspaces $E$ of $\dg X$
and $F$ of $\dg Y$ and a linear $S$-isomorphism $\psi : E\rightarrow F$ which
is a similitude of the quadratic spaces.

(i)
$S/(T \times U)$ embeds as
a diagonal subgroup of
$O(E) \times O(F)$, i.e., a subgroup which meets each direct factor
$O(E)$ and $O(F)$ trivially.

(ii)
The image of $S$ in $O(F)$ is
$A\cap B$ where $A$ is the image of $Stab_{O(Y)}(F)$ in $O(F)$ and
$B$ is the subgroup
obtained  from $\psi $ of
the image $C$ of $Stab_{O(X)}(E)$ in $O(E)$, i.e.,
$B=\psi C \psi^{-1}$.  (A similar description applies to the projection to $O(\dg X)$).

(iii) We use the notations of (i, ii).
Suppose that $E=\dg X$ and $F=\dg Y$.
 If $O(X)$ induces $O(\dg X)$ on $\dg X$ and
$O(Y)$ induces $O(\dg Y)$ on $\dg Y$, then
$S/(T \times U)\cong O(\dg X)\cong O(\dg Y)$.

\end{lem}
\pf Straightforward. \eop

\subsection{About the centralizer of an $EE_8$ dihedral group of order 10 in $O({\mathcal N}(A_4^6)$}

%

Let $\mathcal{N}:=\mathcal{N}(A_4^6)$
be the Niemeier lattice with the root
sysytem $A_4^6$. Then $\mathcal{N}\supset A_4^6$. Since $A_4^*/A_4\cong
\ZZ_5$, $\mathcal{N}/A_4^6< \ZZ_5^6$ is a self-dual  linear code over $\ZZ_5$
and it is generated by
\[
(1,0,1,-1,-1,1) , \quad (1,1, 0,1,-1,-1), \quad (1,-1,1, 0,1,-1)
\]
(see \cite{splag}, \cite{Erokhin}).

\medskip

Let $H$ be the automorphism group of the glue code $\mathcal{C}
=\mathcal{N}/A_4^6$. Then $H< \ZZ_2^6. Sym_6$. Now use $\{\infty, 0,1,2,3,4\}$
to label the 6 coordinates . Then $H$ has shape $2.PGL(2,5)$ and
generated by
\[
(01234),\quad
\varepsilon_2 \varepsilon_3 (\infty\, 0)(23), \quad
\varepsilon_{\infty} (1243),\quad
\varepsilon_\infty \varepsilon_0\varepsilon_1\varepsilon_2 \varepsilon_3 \varepsilon_4,
\]
where $\varepsilon_{i}$ is the mutliplication of $-1$ on the $i$-th coordinate
and a permutation denotes a permutation matrix (cf. \cite{splag},
\cite{Erokhin} and \cite{lps}). Note that a cycle $(i_1i_2 \dots i_\ell)$ acts
on $\{\infty, 0,1,2,3,4\}$ from the left and maps $i_1$ to $i_2$, $i_2$ to  $
i_3$, \dots, $i_\ell$ to $i_1$.

\medskip

Note that $H$ also acts on $\mathcal{N}$ and the isometry group
$O(\mathcal{N})$ is isomorphic to the semidirect product $ W \rtimes H
$, where $W= Weyl(A_4^6)\cong Sym_5^6$ and $H\cong 2.PGL(2,5)$ is the
automorphism group of the glue code (cf. \cite{splag}, \cite{Erokhin}).

\medskip

Let $h$ be a fixed point free order 5 element in $O(A_4)$  and define
\[
M:= span \{ (0,0,\a, \b, -\b, -\a)\mid \a, \b\in A_4\}\cup \{ (0,0,\gamma,
\gamma', -\gamma', -\gamma)\}
\] and
\[
N:= span \{ (0,0,h\a, h\b, -\b, -\a)\mid \a, \b\in A_4\}\cup \{  (0,0,h\gamma,
h\gamma', -\gamma', -\gamma)\},
\]
where $\gamma=\frac{1}5(1,1,1,1,-4)$ and $\gamma'= \frac{1}5(2,2,2,-3,-3)$.
Then $M\cong N \cong EE_8$. Notice that the word $(0,0,1,2,-2,-1)$ is in the
glue code.

\medskip

The corresponding SSD involutions are given by
\[
t_M= (14)(23)
\]
and
\[
t_N= (1, 1, h, h, h^{-1}, h^{-1}) (14)(23).
\]
Thus, $t_Nt_M = (1, 1, h, h, h^{-1}, h^{-1})$ has order $5$  and  the dihedral
group $\la t_M, t_N\ra$ is isomorphic to $Dih_{10}$ and is generated by
\[
\tilde{h}:= (1, 1, h, h, h^{-1}, h^{-1})\quad \text{ and } \quad t_M=(14)(23).
\]


We define $Q:=M+N$, and $R':=ann_{\mathcal N}(Q)$ and $Q':=ann_{\mathcal
N}(R')$.
 Since $t$ and $u$ act on the six coordinate spaces for the glue code as permutations of cycle shape $1^22^2$,
$R'$ contains the sum of two $A_4$ indecomposable orthogonal components of the root sublattice.
 Then $Q'\perp R'$ contains a copy of $Q\perp R$
 and  $\mathcal N$ is an overlattice of $Q\perp R$ as described in \refpp{el3}, the case $e=1$.
 In particular $R'\cong A_4^2$, $\dg {R'} \cong 5^2 \cong \dg{Q'}$ and the root sublattice of $Q'$ is isometric to $A_4^4$.

 If $L$ is an integral lattice, we let $Weyl(L)$ denote the group generated by reflections at the roots of $L$.

Clearly, $C_{Weyl(Q'\perp R')}(\tilde{h}) \cong Weyl(A_4^2) \times 5^4$.
 Also,
$C_{O(\mathcal N)}(\tilde{h})$ is contained in the subgroup of index $6 \choose
2$ of $O(\mathcal N)$ which fixes the two orthogonal direct summands of
$Q'\perp R'$ which are in $R'$.   The group $C_{O(\mathcal N)}(\tilde{h})$
therefore has the property that its  image in $O(\mathcal N)/Weyl(Q'\perp
R')\cong 2.PGL(2,5)$ has order $240/{6\choose 2}=16$.  It follows that
$C_{O(\mathcal N)}(\tilde{h})^{(\infty )}$, the terminal member of the derived
series of $C_{O(\mathcal N)}(\tilde{h})$, is just $Weyl(R')'\cong Alt_5 \times
Alt_5$. This is contained in $C_{O(\mathcal N)}(D)$ and in fact equals
$C_{O(\mathcal N)}(D)^{(\infty )}$.

We have proved:

\begin{prop}\labtt{centdniema4to6}
We let $\mathcal N := \mathcal N (A_4^6)$ and let $D$ be a dihedral group of order 10 generated by a pair of $EE_8$-involutions in $O(\mathcal N )$ which give permutations of cycle shape $1^22^2$ on the six coordinate space for the glue code (equivalently, on the six indecomposable orthogonal components of the root sublattice).  Then
$C_{O(\mathcal N)}(D)^{(\infty )} = C_{O(\mathcal N)}(O_{2'}D)^{(\infty )} \cong Alt_5 \times Alt_5$.
\end{prop}

\section{The finite orthogonal groups}

We give formulas for the orders of the orthogonal groups over finite fields.

\begin{thm}\labtt{orderog}
We let $q$ be an odd prime power, $\vep = \pm$.

(i)
$|SO^{\vep }(2n,q)|=q^{n(n-1)} (q^n-\vep ) \prod _{i=1}^{n-1}(q^{2i}-1)$

$|SO(2n+1,q)|=q^{n^2} \prod _{i=1}^{n}(q^{2i}-1)$

\smallskip

(ii) The orders of the corresponding orthogonal groups are twice the above numbers.
\end{thm}
\pf
(i)
This is taken from \cite{artin}, pages 126-147.  (ii) In characteristic not 2, an  orthogonal group contains the special orthogonal group with index 2.
\eop

\begin{coro}\labtt{stabmaxtotsing}
The stabilizer in $O^+(2n,q)$ of a maximal totally singular subspace has the form
$q^{n \choose 2}{:}GL(n,q)$.  The number of such subspaces is
$2\prod_{i=1}^{n-1}(q^i+1)$.
\end{coro}
\pf
The form of the stabilizer follows by considering a direct sum decomposition of the space by a pair of totally singular subspaces.  The second statement follows from Witt's transitivity result and the formulas of \refpp{orderog}.
\eop

\begin{coro}\labtt{fourdimog }
Suppose that $q$ is odd and that  $\vep = +$.  The orthogonal group has structure
$(SL(2,q)\circ SL(2,q)).(2 \times 2)$.  The normalizer of either $SL(2,q)$ central factor is the special orthogonal group.
\end{coro}
\pf
By order considerations \refpp{orderog}, it suffices to show that
the orthogonal group contains such a subgroup.

For $i=1,2$, let $U_i$ be $\FF_q^2$ with a nonsingular alternating form.  The
group of similitudes is $GL(U_i)$.   The wreath product $(GL(U_1)\times
GL(U_2))\la t \ra$, where $t$ is the usual involution which switches arguments of
tensors,  acts on $V:=U_1\otimes U_2$, as simitudes of the nonsingular
quadratic form obtained by tensor product from the forms on the $U_i$.   Call
this homomorphism $\varphi$. The kernel $K$ of $\varphi$ is a central subgroup
isomorphic to $\ZZ_{q-1}$. The subgroup $H$ of $GL(U_1)\times GL(U_2)$
preserving the form is normal and gives a quotient isomorphic to $\ZZ_{q-1}$.
Obviously, $K\le H$.

Let $D_i$ be the diagonal subgroup of $GL(U_i)$, $i=1,2$.
Write $h_i(a,b)$ for the diagonal matrix
$\begin{pmatrix}a&0 \cr 0&b \end{pmatrix}$ of $D_i$.

Define $D=D_1 \times D_2$.
Then $D\cap H = \{ (h_1(a,b), h_2(c,d)) \mid abcd=1 \}$
so $H$ is just the kernel of the map
$\psi : GL(U_1)\times GL(U_2) \rightarrow \FF_q^{\times}$ which takes
$(g_1, g_2)$ to $det(g_1)det(g_2)$.

The intersection $Z$ of
$Z(GL(U_1)) \times Z(GL(U_2))$
 with $H$ is
$Z= \{ h_1(a,a)h_2(b,b) \mid (ab)^2=1, \text{i.e.,} ab=\pm1 \}$
and $K = \{ h_1(a,a)h_2(b,b) \mid ab=1 \}$.

Fix a nonsquare $n\in \FF_q^{\times}$.
Observe that, given $a, b$, there is $e$ so that $abe^2 \in \{1, n\}$.
It follows that if $(h_1(a,b), h_2(c,d))  \in D\cap H$,
it is congruent modulo $K$ to
$(h_1(1,1), h_2(1,1))$ or
$(h_1(n,1), h_2(n^{-1},1))$.
Therefore,
$Im(\varphi )$ contains $\varphi (SL(U_1)\times SL(U_2))$ with index 2.  By order considerations, $Im(\varphi )=SO(V)$.

Finally, note that the action of $t$ has determinant $-1$ and that it interchanges the groups
$\varphi (SL(U_1))$ and $\varphi (SL(U_2))$ under conjugation.
\eop

\begin{coro}\labtt{fourdimog2}  Let $q$ be an odd prime power.
A subgroup $S$ of index 2 in $O^+(4,q)$ contains $Z:=Z(O^+(4,q))=\{\pm 1\}$ and
$S/Z$ is isomorphic to
either $PSL(2,q)\wr 2$ or
the even subgroup of $PGL(2,q)\times PGL(2,q)$, i.e., the subgroup of index 2 which intersects each direct factor in its unique subgroup of index 2.
\end{coro}
\pf  The containment $Z\le S$ follows since $Z$ is in the commutator subgroup of $O^+(4,q)$.

We use
the fact that
$PSO^+(4,q)$ embeds as a subgroup $T$ of index 2
in $U:=PGL(2,q)\wr 2$.  We have
$U''\cong PSL(2,q)\times PSL(2,q)$ and $U/U''\cong Dih_8$, in which $T/U''$ is a four-group.  In $Dih_8$, there are just two conjugacy classes of four-groups.
Our group $T$ is not the
subgroup
$PGL(2,q)\times PGL(2,q)$.  Furthermore, the two subgroups of order 2 in
$T/U''$ which are not central in $U/U''$ are conjugate.
Therefore, the isomorphism types for $S$ are limited to two possibilities.
Now use  \refpp{fourdimog }.
\eop

\begin{lem}\labtt{threedimog}
Let $q$ be an odd prime power.  Then $O(3,q)\cong 2\times PGL(2,q)$.
\end{lem}
\pf
Since the dimension is odd, $O(3,q)=\la -1 \ra \times SO(3,q)$. It suffices to show that
$SO(3,q)\cong PGL(2,q)$.

We take the natural actions of $GL(2,q)$ on
the space of $2\times 2$ matrices over $\FF_q$ and observe that it leaves invariant the trace 0 matrices for which the form $A, B\mapsto Tr(A, B)$ is nonsingular and symmetric.  The kernel of this action is the group of scalar matrices.  Consequently, $PGL(2,q)$ embeds in $SO(3,q)$.  The order formula \refpp{orderog} shows that this embedding is onto.
\eop

\begin{lem}\labtt{q6}
Let $q$ be an odd prime power.
There are five orbits of $O^+(4,q)$ on 2-spaces, and these orbits are distinguished by these properties of their members (in the notation $Orbit(a,\cdots)$, the parameter $a$ means the dimension of the radical of the 2-space):

Orbit(2): totally singular;

Orbit(1,s), Orbit(1,n):  nullity 1;  two kinds, according to whether norms are squares or nonsquares

Orbit(0,1): nonsingular, maximal Witt index;

Orbit(0,0): nonsingular, nonmaximal Witt index.
\end{lem}
\pf
This is a consequence of Witt's theorem and the theory of nonsingular quadratic spaces.
\eop

\begin{lem}\labttr{j2tensorj2}
Let $p$ be an odd prime and let $\la g \ra$ be a group of order $p$.  Suppose that for $i=1,2$,  $U_i$ is a 2-dimensional $\FF_p\la g \ra$-module on which  $g$ acts with a single degree 2 indecomposable Jordan block.
The action of $g$ on $U_1 \otimes U_2$ has Jordan canonical form the sum of a degree 1 and degree 3 indecomposable Jordan block.
\end{lem}
\pf
For $U_i$, choose a basis $e_i, f_i$ so that $g$ fixes $f_i$ and $ge_i=e_i+f_i$.
Consider the   $\FF_p\la g \ra$-submodules $span\{e_1\otimes e_2, f_1\otimes e_2+e_1\otimes f_2, f_1\otimes f_2\}$ and
$span\{ f_1\otimes e_2 - e_1\otimes f_2\}$.
\eop

\begin{rem}\labttr{jmtensorjn}
Tensor products of indecomposables for cyclic $p$-groups in characteristic $p$ were studied in \cite{srinivasin,lindsey}.
\end{rem}

\begin{lem}\labttr{r13.5} In $O^+(4,q)$, for $q$ a power of the odd prime $p$,
the elements of order $p$ in each $SL(2,q)$-component of $O^+(4,q)$ have Jordan canonical form a direct sum of two indecomposable degree 2 blocks.  The other elements of order $p$ in $O^+(4,q)$ have Jordan canonical form a sum of indecomposable blocks of degrees 1 and 3.
\end{lem}
\pf
\refpp{j2tensorj2}.  \eop

\begin{lem}\labtt{q6.5}
We use the hypotheses and notations of \refpp{q6}.
Suppose that $v_1, v_2, v_3, v_4$ is a basis of the nonsingular quadratic space
$V :=\FF_q^4$ of maximal Witt index which has Gram matrix
$$\begin{pmatrix}
0&0&0&1\cr
0&a&0&0\cr
0&0&b&0\cr
1&0&0&0\cr
\end{pmatrix}, a \ne 0, b\ne 0.$$
Let $W:=span\{v_1, v_2\} \in Orbit(1,s)\cup Orbit(1,n)$.
The stabilizer in $O^+(4,q)$ of $W$ has the form
$S=UH$, where $U$ is a normal subgroup of order $q^2$ which acts trivially on the factors of the chain
$0 < W < ann(v_1) <V$ and where $H$ is diagonal with respect to the above basis, acting as matrices of the form
$$\begin{pmatrix}
c&0&0&0\cr
0&d&0&0\cr
0&0&e&0\cr
0&0&0&c^{-1}\cr
\end{pmatrix}, c \in \FF_q^{\times}, d, e \in \{\pm 1\}.$$
In particular, $H\cong \ZZ_{q-1} \times 2 \times 2$.  The order of $S$ is $2^2(q-1)q^2$.  The kernel of the action of $S$ on $W$ is a subgroup of order $q$.
\end{lem}
\pf It is clear from the proof of \refpp{j2tensorj2} and \refpp{r13.5}  that $S$
contains a maximal unipotent group, $U$, and that $U$ acts trivially on the
factors of such a chain. Furthermore, any orthogonal transformation which acts
trivially on the factors of such a chain is unipotent so is in $U$.

We have $S \le Y$, the stabilizer of $\FF_q v_1$.  The group $Y$ splits over $U$
by coprimeness.  A complement  $H$ to $U$ in $Y$  has a direct product
decomposition $H=H_1\times H_2$ where $H_1$ acts trivially on
$W':=span\{v_2, v_3\}$ and faithfully on $W'':=span\{v_1, v_4\}$ and such that
$H_2$ acts trivially on $W''$ and faithfully on $W'$.  We have $H_1 \cong
\FF_q^{\times}$ and $H_2\cong O^+(2,q)\cong \FF_q^{\times}:2$, a dihedral
group.   Thus $S$ contains $UH_1$ and $S\cap H_2$ is just the subgroup fixing a
nonsingular 1-space in the natural 2-dimensional orthogonal representation of
this group.  Therefore, $S\cap H_2\cong 2\times 2$.  Our matrix representation
of $S$ follows since $H_2$ fixes a complement to $W'$ in $V$ and we may as
well assume that it is $W''$ by replacing $H$ with a $U$-conjugate.   The final
two statements are easy. \eop

\section{Finite subgroups of $E_8(\CC )$}

\begin{lem}\labtt{r4} A sublattice in $E_8$ of index 5 contains roots.
\end{lem}
\pf
Such a sublattice corresponds to a cyclic group of order 5 in $E_8 (\CC )$.  The list of elements of order 5 shows that every one fixes root spaces. See \cite{cg,grelab}.  \eop

%

%


\begin{lem}\labtt{r4.7}
Let $h\in O(E_8)$ be an isometry of order 5 without eigenvalue 1.  Then $(h-1)E_8$ is rootless.
\end{lem}
\pf Let $E_8$ contain a sublattice $M_1\perp M_2$, where $M_1 \cong M_2 \cong
A_4$. Then we may take $h=h_1h_2$, where $h_i \in Weyl(M_i)$ and $|h_i|=5$.
Then $E_8$ contains $M_1 \perp M_2$ with index 5 and is represented by a
nontrivial gluing. The action of $(h-1)$ carries $E_8$ to a sublattice of
$M_1\perp M_2$ which intersects $M_i$ in the rootless sublattice $(h_i-1)M_i$,
for $i=1,2$ \refpp{r3}. \eop

\begin{lem}\labtt{r5}  We use the notation of \cite{cg} for conjugacy classes in $E_8(\CC )$.
(i) $E_8(\CC )$ has one conjugacy class of subgroups $F$ such that $F\cong 5^2$ and $dim(C(F))=8$, i.e., $C(F)^0$  is a torus.  In fact, $C(F)\cong \TT\cdot 5$.  If $g\in C(F)\setminus C(F)^0$, the action of $g$  on the root lattice modulo 5 has Jordan canonical form the sum of two degree 4 Jordan blocks.

(ii) There is one $E_8(\CC )$-orbit on ordered pairs $(x,y)$ of class  $5C$-elements which generate an elementary abelian group of order 25 whose centralizer is 8-dimensional.

(iii) In the lattice $Y:=E_8$, there is one equivalence class of rootless sublattices $X$ so that $Y/X \cong 5^2$.   We have $Stab_{O(Y)}(X)\cong 5{:}GL(2,5)\cong (5\times SL(2,5)){:}4$, a group with trivial center,
and the normal subgroup  of order 5 is the kernel of the action on $Y/X$.

(iv) An $X$ as in (iii) contains a sublattice isometric to $A_4(1) \perp A_4(1)$.  Furthermore, an isometry $g$ of order 5 in $O(X)\cap O(Y)$ acts nontrivially on $\dg X$, with Jordan canonical form a sum of two degree 2 Jordan blocks.  The Smith invariants of $g-1$ on
$X$ are $(1,1,1,1,1,1,5,5)$.
\end{lem}
\pf
The arguments use techniques from the theory of finite subgroups of Lie groups.  See \cite{gr99,gr02} for references and examples.  Let $Tr$ be the trace function for a linear transformation on the adjoint module for the Lie group of type $E_8$.

(i)  and (ii)  A survey of the elements of finite order (e.g., \cite{gr99,gr02}) shows
that for a cyclic group $H$ of order 5 in $E_8(\CC )$, $\sum_{g\in H} Tr(g) \ge
240$, with equality occurring only at the case where nonidentity elements of $H$
are in class $5C$ (recall that this sum is the $|H|$ times the dimension of the
fixed point subalgebra). Consequently, if $E$ is any elementary abelian group of
order 25, $\sum_{g\in E} Tr(g) \ge  200$, with equality holding only when all
nonidentity elements of $E$ are in the class $5C$.

(iii) This follows from the proof of (ii), which shows that if $x$ has type $5C$, it
corresponds to a sublattice $Z$ of $Y$ of type $A_4^2$.  A rootless index 5
sublattice of this may be obtained as follows.  Let $Z_1$ and $Z_2$ be the pair of
$A_4$ direct summands. For $i=1,2$, let $g_i$ be an element of order 5 from the
Weyl group of $Z_i$. Any sublattice between $(g_1 -1)Z_1\perp (g_2-1)Z_2$ and
$Z_1\perp Z_2$ is stable under the action of $g_1g_2$.  In particular, a
sublattice $J$ of index 5 in $Z_1\perp Z_2$ which meets $Z_i$ in $(g_i-1)Z_i$, for
$i=1,2$, is rootless and is left invariant by $g_1g_2$.
There are four such lattices $J$, and this set of four lattices forms an orbit under
$W(Z_1)\times W(Z_2)$, where $W(Z_i)\cong Sym_5$ denotes the Weyl group.
Since $W(Z_1)\times W(Z_2) \le O(Y)$, it follows that
$g_1g_2$ is conjugate to all of its nontrivial powers in $O(J)\cap O(Y)$.

Such a sublattice $J$ corresponds to an elementary abelian group of order 25 which fixes no root spaces, hence whose connected centralizer is a maximal torus.
By (ii), its normalizer in $E_8(\CC)$ induces on $Y/X$ the group $GL(2,5)$ and so the kernel of the action has shape $\TT.5$, whose component group is cyclic of order 5.
Now apply (ii) to $J$ and the action of $g_1g_2$ on $J$.

(iv) We use the notation of (iii).  Let $P_i$ be the orthogonal projection to $\QQ
Z_i$, $i=1,2$.

Since $g_i$ has minimum polynomial $(x^5-1)/(x-1)\in \QQ [x]$ on $Z_i$,
it has minimum polynomial $(x-1)^4 \in \FF_5 [x]$ on $A/5A$, where $A$ is any nonzero
$g_i$-submodule of $Z_i$ \refpp{r4}(iii).

Now set $g:=g_1g_2$.
Then the minimum polynomial of $g$ on any nonzero invariant lattice
$B$ in $\QQ R$ is $(x^5-1)/(x-1)\in \QQ [x]$.   By reasoning as above, we get that the minimum polynomial of $g$ on $B/5B$ is $(x-1)^4\in \FF_5 [x]$.

The central involution of $Stab_{O(Y)}(X)\cong 5\times GL(2,5)$ acts as $-1$ so
acts on $\dg B$ with irreducible constituents of even degree. Since, by minimal
polynomial considerations, $g$ acts nontrivially on $\dg B$, it follows that the
commutator space of $g$ on $\dg B$  has dimension 2 and so $g$ acts with
Jordan canonical form the sum of two indecomposable degree 2 blocks.

By
\refpp{r4.5}, we see that $(g-1)^4X=5X$, which implies that $5^2$ does not divide any Smith invariants.  Since $det(g-1)=25$, the result follows.
 \eop

\begin{lem}\labtt{r10.5}  If $g$ is an element of $E_8(\CC )$ of finite order and $dim(C(g))=8$, then $|g|\ge 30$.
\end{lem}
\pf  See \cite{cg,grelab,gr99,gr02} for background on finite subgroups of Lie groups.
An element of finite order is conjugate to one which corresponds to a labeling of the extended $E_8$-diagram \cite{kac}.  The condition $dim(C(g))=8$ implies that each label is positive, whence $|g|\ge 30$.
\eop

\begin{lem}\labtt{qine8e8.2}
Let $U_1\perp U_2$ be an orthogonal direct sum of $E_8$-lattices.
Suppose that  $Q \le U_1\perp U_2$ and $Q$ is a rootless lattice of index $5^2$.
Then
 $Q \le U_1\perp U_2$ corresponds to an elementary abelian
group $E\cong 5^2$ in the Lie group $G=E_8(\CC )\times E_8(\CC )$ whose
connected centralizer in $G$ is a torus.  Such a group $E$ is unique in $G$ up to
conjugacy.
\end{lem}
\pf Let $G_i, i=1,2$ be the two direct factors of $G$ and let $E_i$ be the
projection of $E$ to $G_i$.  Then the connected centralizer of $E_i$ in $G_i$ is a
torus, $i=1,2$.  By \refpp{r8}, $E_i$ is unique up to conjugacy in $G_i$ and its
normalizer $N_i$ induces $Aut(E_i)\cong GL(2,5)$ on $E_i$, for $i=1, 2$. The
action of $N_1\times N_2$ on $E_1 \times E_2$ has one orbit on its subgroups
of order $5^2$ which meet each of $E_1$ and $E_2$ trivially. The lemma follows.
\eop

\section{Finite subgroups of $O(n, \mathbb{F})$ which centralize root groups}

Let $\mathbb{F}$ be an algebraically closed field.

\begin{lem}\labtt{centrootgroup}
Let $A$ be an odd order finite abelian subgroup of $GL(V)$  where $V$ is a $2m$
dimensional vector space over $\mathbb F$, an algebraically closed field of characteristic
not dividing $|A|$. We suppose that $A$ leaves invariant a split, nondegenerate
quadratic form. Denote by $O(V)$ the isometry group of this form.

(i) If there exists a nontrivial linear character of $A$ which occurs with
multiplicity at least 2, then $C(A)$ contains a root group.

(ii) If the trivial character occurs with multiplicity at least 4, then $C(A)$ contains
a root group.
\end{lem}
\pf  (i) For a linear character $\mu$ of $A$, let $V(\mu)$ be the eigenspace for
$\mu$ in $V$.

Let $\l \ne 1$ be a character which occurs with multiplicity $k\ge 2$. Then $V(\l
)$ and $V(\l^{-1})$ each have dimension $k$ and intersect trivially.  Furthermore,  since $|A|$ is odd,
the spaces $V(\l )$ and $V(\l^{-1})$ are each totally singular and the associated
bilinear form pairs them nonsingularly.

Let $W$ be an $m$-dimensional totally singular subspace which is invariant
under $A$ and which contains $V(\l)$. Let $W'$ be an $m$-dimensional totally
singular subspace which is invariant under $A$ and which contains $V(\l^{-1})$.
We may arrange for $W\cap W'=0$. The subgroup $H$ of $O(V)$ which leaves
both $W$ and $W'$ invariant acts faithfully on each of $W$ and $W'$ as
$GL(W), GL(W')$, respectively.

Since $k\ge 2$, $A$ centralizes a natural $GL(k,\mathbb F)$-subgroup of $H$ and so
centralizes a root group of $O(V)$.

(ii) This is immediate since a natural 4-dimensional orthogonal group contained in
$O(V)$ contains a root group. \eop


\begin{coro}\labttr{elabcentrootgroup}
Let $p$ be an odd prime and let $A$ be an elementary abelian  group of order
$p^2$ in $O(2m,\mathbb F)$, where $\mathbb F$ is an algebraically closed field of characteristic not dividing $|A|$.  If $p^2< 4m-3$, then $C(A)$ contains a root group.
\end{coro}
\pf By \refpp{centrootgroup}, we may assume that the fixed point space has
dimension at most 2.  We may partition the nontrivial linear characters into
inverse pairs.  If one such pair occurs with multiplicity at least 2, we are done by
\refpp{centrootgroup}.  Denying this, we get $\frac {p^2-1}2 \ge 2m-2$, or
$p^2\ge 1+4m-4=4m-3$. \eop

\begin{coro}\labttr{qinhs16}
Let $A$ be an abelian group of order 25 in $O(16,F)$, where
$\mathbb F$ is an algebraically closed field of characteristic not 5.
Then $A$ centralizes a root
group.
\end{coro}
\pf We note that $25<32-3$ then use \refpp{elabcentrootgroup}. \eop


\end{document}